\documentclass{elsart3-1}

\usepackage{amssymb}
\usepackage{amsmath}
\usepackage{mathrsfs}
\usepackage{paralist}
\usepackage{dsfont}
\usepackage{cite}
\usepackage{setspace}

\usepackage[english,francais]{babel}

\newcommand{\calH}{{\mathcal{H}}}
\newcommand{\calL}{{\mathcal{L}}}
\newcommand{\eps}{\varepsilon}
\newcommand{\W}{\mathbf{W}}
\newcommand{\m}{\mathbf{m}}
\newcommand{\N}{{\mathbb{N}}}
\newcommand{\R}{{\mathbb{R}}}

\newcommand{\dd}{\,\mathrm{d}}

\newcommand{\dff}{\mathrm{D}}

\newcommand{\dist}{\mathrm{\mathbf{dist}}}
\newcommand{\prb}{\mathscr{P}_2}
\newcommand{\prbb}{\mathscr{P}}

\newcommand{\theX}{\mathbf{X}}

\newcommand{\au}{\mathcal{U}}
\newcommand{\av}{\mathcal{V}}
\newcommand{\ent}{\mathcal{E}}

\newcommand{\equ}{u^\infty}
\newcommand{\eqv}{v^\infty}
\newcommand{\matr}[1]{\mathbf{#1}}

\newtheorem{theorem}{Theorem}[section]

\newtheorem{e-proposition}[theorem]{Proposition}

\newtheorem{e-definition}[theorem]{Definition\rm}

\renewcommand{\theequation}{\arabic{equation}}
\def\og{\leavevmode\raise.3ex\hbox{$\scriptscriptstyle\langle\!\langle$~}}
\def\fg{\leavevmode\raise.3ex\hbox{~$\!\scriptscriptstyle\,\rangle\!\rangle$}}

\allowdisplaybreaks

\journal{the Acad\'emie des sciences}
\begin{document}
% place in the next line the header (rubrique) chosen for your article,
% if you know it (you can also have 2, format : Header1/Header2
\centerline{Partial Differential Equations}
\begin{frontmatter}


\selectlanguage{english}
\title{A note on the variational analysis of the parabolic-parabolic Keller-Segel system in one spatial dimension}
\date{August 14, 2014}

\selectlanguage{english}
\author[JZ]{Jonathan Zinsl},
\ead{zinsl@ma.tum.de}

\address[JZ]{Zentrum f\"ur Mathematik, Technische Universit\"at M\"unchen, 85747 Garching, Germany}

% If you know the dates of reception, and acceptation you can put them now;
%  idem the name of the person presenting the Note

\medskip
\begin{center}
{\small Received *****; accepted after revision +++++\\
Presented by xxxxx}
\end{center}

\begin{abstract}
\selectlanguage{english}
% Text of abstract in English
We prove the existence of global-in-time weak solutions to a version of the parabolic-parabolic Keller-Segel system in one spatial dimension. If the coupling of the system is suitably weak, we prove convergence of those solutions to the unique equilibrium with an exponential rate. Our proofs are based on an underlying gradient flow structure with respect to a mixed Wasserstein-$L^2$ distance.

\vskip 0.5\baselineskip

\selectlanguage{francais}
% Text of abstract in French
\noindent{\bf R\'esum\'e} \vskip 0.5\baselineskip \noindent
{\bf Une note sur l'analyse variationelle du syst\`eme de Keller-Segel parabolique-parabolique \`a une dimension spatiale. }
Nous prouvons l'existence de solutions faibles globales en temps d'une variante du syst\`eme de Keller-Segel parabolique-parabolique \`a une dimension spatiale. Si le couplage du syst\`eme est assez faible, nous prouvons la convergence de ces solutions vers l'equilibre univoque \`a une vitesse exponentielle. Nos preuves se reposent sur une structure de flux de gradient dans l'espace produit des espaces Wasserstein et $L^2$.

\end{abstract}
\end{frontmatter}

\selectlanguage{english}

\section{Introduction and main results}

We consider the following version of the Keller-Segel model for chemotaxis in one spatial dimension:
\begin{align}
u_t(t,x)&=(u_x(t,x)+u(t,x)W_x(x)-\chi u(t,x)v_x(t,x))_x,\label{eq:pde_u}\\
v_t(t,x)&=v_{xx}(t,x)-\kappa v(t,x)+\chi u(t,x),\label{eq:pde_v}
\end{align}
where $t>0$ and $x\in \R$, and the sought solution $(u,v)$ is subject to the initial conditions
\begin{align}\label{eq:ic}
u(0,x)=u^0(x),\qquad v(0,x)=v^0(x).
\end{align}
We require that $\kappa\ge 0$, $\chi\in\R$ and the \emph{confinement potential} $W\in C^2(\R)$ to be bounded from below and to grow at most quadratically, i.e. $\underline{W}\le W(x)\le Ax^2+B$ for all $x\in \R$ and some $\underline{W},A,B\in \R$, and to have bounded second derivative $W_{xx}$.

It is known that \eqref{eq:pde_u}\&\eqref{eq:pde_v} possesses a variational structure since it can formally be written as a gradient flow of the (non-convex) entropy functional $\calH:\,\theX\to\R_{\infty}$ (see formula \eqref{eq:ent} below) with respect to the compound distance 
%\begin{align*}
$\dist((u,v),(u',v')):=\sqrt{\W_2^2(u,u')+\|v-v'\|_{L^2}^2}$
%\end{align*}
on the space $\theX:=\prb(\R)\times L^2(\R)$, where $(\prb(\R),\W_2)$ is the space of (absolutely continuous) probability measures -- or their densities, respectively -- on $\R$ with finite second moment $\m_2$, endowed with the $L^2$-Wasserstein distance $\W_2$. The entropy $\calH$ is defined as
\begin{align}\label{eq:ent}
\begin{split}
\calH(u,v)&=\begin{cases}\int_\R\left[u\log u+uW+\frac12 v_x^2+\frac{\kappa}{2}v^2-\chi uv\right]\dd x,&\text{if }\int_\R u\log u\dd x<\infty\text{ and }v\in H^1(\R),\\
+\infty,&\text{otherwise.}\end{cases}
\end{split}
\end{align}

In this note, we sketch another possible application of the method in \cite{zinsl2014} to prove the existence of weak solutions to \eqref{eq:pde_u}\&\eqref{eq:pde_v} and to analyse their long-time behaviour. There, it has been shown that in the case of a porous-medium-type diffusion for $u$ on $\R^3$, a global-in-time weak solution always exists, and that it converges exponentially fast to the unique equilibrium if the coupling of the system is suitably weak. In the one-dimensional setting at hand, some parts of the proofs simplify compared to those in \cite{zinsl2014} due to a gain in regularity. In contrast to that, the case of linear diffusion causes the difficulty of a missing time-uniform \emph{a priori} estimate for $u$ in $L^m(\R)$ for some $m>1$.

The cornerstone of our variational analysis is the so-called \emph{minimizing movement scheme} (see e.g. \cite{savare2008, jko1998}) for the construction of an approximate time-discrete gradient flow w.r.t. the distance $\dist$:

For each step size $\tau>0$,
let $(u_\tau^0,v_\tau^0):=(u^0,v^0)$, and then define inductively for each $n\in\N$:
\begin{align}
  \label{eq:jko}
  (u_\tau^n,v_\tau^n) \in \operatorname*{argmin}_{(u,v)\in\prb(\R)\times L^2(\R)} \Big(\frac1{2\tau}\dist\big((u,v),(u_\tau^{n-1},v_\tau^{n-1})\big)^2 + \calH(u,v)\Big).
\end{align}

Further, introduce the piecewise constant interpolation $(u_\tau,v_\tau):\R_+\to\prb(\R)\times L^2(\R)$ by
\begin{align}
  \label{eq:interpol}
  u_\tau(t) = u_\tau^n,\quad v_\tau(t) = v_\tau^n \quad \text{for all $t\in((n-1)\tau,n\tau]$}.
\end{align}

This hybrid variational principle has been exploited previously for Keller-Segel-type systems \cite{blanchet2014, blanchet2013, mimura2012, zinsl2013} in higher spatial dimensions and also in other applications, e.g. \cite{kinderlehrer2009, laurencot2011}. For the vast literature on the behaviour of the Keller-Segel system and its variants, we refer to the review articles by Horstmann \cite{horstmann2003} and Blanchet \cite{blanchet2013r} and emphasize that the one-dimensional model on bounded spatial domains has been explicitly investigated by Osaki and Yagi \cite{osaki2001} and Hillen and Potapov \cite{hillen2004}, leading to similar results as proven here.

We obtain the following statement on the existence of global-in-time weak solutions:
\begin{theorem}[Existence]\label{thm:exist}
Assume that $\chi$, $\kappa$ and $W$ are as mentioned above and that the initial condition satisfies $u^0\in \prb(\R)$, $\int_\R u^0\log u^0\dd x<\infty$ and $v^0\in H^1(\R)$. Define, for each $\tau>0$, a discrete solution by means of \eqref{eq:jko}\&\eqref{eq:interpol}. Then, there exists a vanishing sequence $\tau_k\searrow 0~(k\to\infty)$ such that $(u_{\tau_k},v_{\tau_k})$ converges to a weak solution $(u,v)$ to \eqref{eq:pde_u}--\eqref{eq:ic} in the sense that \eqref{eq:pde_u} holds in the sense of distributions, whereas \eqref{eq:pde_v} and \eqref{eq:ic} hold almost everywhere. One has for all $T>0$:
\begin{align*}
	u_{\tau_k}&\rightharpoonup u \text{ narrowly in the space of probability measures $\prbb(\R)$, pointwise with respect to $t\in[0,T]$},\\
	v_{\tau_k}&\to v \text{ in $L^2(\R)$, uniformly with respect to $t\in[0,T]$},\\
u&\in C^{1/2}([0,T];(\prb(\R),\W_2))\cap L^1([0,T];L^\infty(\R))\cap L^2([0,T];L^2(\R)),\\
	\sqrt{u}&\in L^2([0,T];H^1(\R)),\quad	u\log u\in L^\infty([0,T];L^1(\R)),\\
	v&\in C^0([0,T]\times \R)\cap H^1([0,T];L^2(\R))\cap L^\infty([0,T];H^1(\R))\cap L^2([0,T];H^2(\R)).
  \end{align*}
\end{theorem}

In particular, for fixed $t>0$, $u(t,\cdot)$ is nonnegative, continuous and bounded. The second component $v$ is bounded and continuous in \emph{both} variables, whereas its nonnegativity could also be obtained starting with a nonnegative initial condition.

Our result on the long-time behaviour of the weak solution from Theorem \ref{thm:exist} reads as follows:
\begin{theorem}[Convergence to equilibrium]\label{thm:conv}
Assume in addition to the hypotheses of Theorem \ref{thm:exist} that $W$ is $\lambda_0$-convex for some $\lambda_0>0$ and that $\kappa>0$ is strictly positive. There exist $\bar\eps>0$, $C>0$ and $L>0$ such that for all $\chi=\eps\in (0,\bar\eps)$, the following statements hold:
\begin{enumerate}[(a)]
\item The system \eqref{eq:pde_u}\&\eqref{eq:pde_v} possesses a unique stationary state $(\equ,\eqv)\in (\prb\cap L^\infty)(\R)\times H^2(\R)$ satisfying 
\begin{align*}
\equ&=U_\eps\exp(-W+\eps\eqv),\qquad\text{$U_\eps>0$ such that $\|\equ\|_{L^1}=1$,}\\
\eqv_{xx}&=\kappa \eqv-\eps \equ.
\end{align*}
\item One has $\Lambda_\eps:=\min(\kappa,\lambda_0)-\eps L>0$ and for all $t\ge 0$, the weak solution $(u,v)$ to \eqref{eq:pde_u}--\eqref{eq:ic} from Theorem \ref{thm:exist} admits the estimate
\begin{align}
\begin{split}
\|u(t)-\equ\|_{L^1}&+\W_2(u(t),\equ)+\sup_{x\in\R}|v(t)-\eqv|+\|v(t)-\eqv\|_{H^1}\\&\le C(\calH(u^0,v^0)-\calH(\equ,\eqv))^{1/2}e^{-\Lambda_\eps t},
\end{split}
\label{eq:expconv}
\end{align}
i.e. $(u(t),v(t))$ converges exponentially fast with rate $\Lambda_\eps$ to the equilibrium $(\equ,\eqv)$ as $t\to\infty$.
\end{enumerate}
\end{theorem}

\section{Sketch of proof for Theorem \ref{thm:exist}}

The crucial step in the proof of Theorem \ref{thm:exist} is to verify that the discrete solution $(u_\tau,v_\tau)$ is well-defined and regular enough to allow for passage to the continuous-time limit $\tau\searrow 0$ in a strong sense. Once obtained, we can proceed as in \cite{zinsl2013, zinsl2014} establishing an approximate weak formulation which turns into the weak formulation of the time-continuous equation as $\tau\searrow 0$. We prove the following
\begin{e-proposition}[Minimizing movement]\label{prop:minmov}
For each $\tau>0$ and $(\tilde u,\tilde v)\in \theX$, the functional 
$$\calH_\tau(\cdot|\tilde u,\tilde v):=\frac1{2\tau}\dist^2(\cdot,(\tilde u,\tilde v))+\calH$$
possesses a minimizer $(u,v)\in \prb(\R)\times H^1(\R)$ with $\int_\R u\log u\dd x< \infty$. Moreover, there exist constants $K_0,K_1,K_2>0$ such that
\begin{align}\label{eq:addreg}
\tau\|(\sqrt{u})_x\|_{L^2}^2+\tau \|v_{xx}\|_{L^2}^2&\le K_0\int_\R(u\log u-\tilde u\log\tilde u)\dd x+K_1(\|v\|_{H^1}^2-\|\tilde v\|_{H^1}^2)+K_2\tau (\|v\|_{H^1}^2+1).
\end{align}
In particular, if in addition $\tilde v\in H^1(\R)$ and $\int_\R \tilde u\log \tilde u\dd x<\infty$, then $v\in H^2(\R)$, $\sqrt{u}\in H^1(\R)$ and $u\in L^\infty(\R)$.
\end{e-proposition}

\noindent\textit{Proof.} 
First, in one spatial dimension, there exists $C_0>0$ such that $\|v\|_{L^\infty}\le \|v\|_{C^{0,\frac12}}\le C_0\|v\|_{H^1}$. Moreover, for some $C_1>0$, one has
\begin{align*}
\int_\R u\log u\dd x&\ge -C_1(\m_2(u)+1)^{1/2}.
\end{align*}
From this, we easily see that for all $(u,v)\in\prb(\R)\times H^1(\R)$ with $\int_\R u\log u\dd x<\infty$, we have
\begin{align*}
\int_\R u\log u\dd x+\underline{W}+\frac12\|v_x\|_{L^2}^2-|\chi|C_0\|v\|_{H^1}\le \calH(u,v)<\infty.
\end{align*}
Using the triangle inequality for $\dist$ and Young's inequality, we deduce coercivity of $\calH_\tau(\cdot|\tilde u,\tilde v)$:
\begin{align*}
\calH_\tau(u,v|\tilde u,\tilde v)&\ge \frac14 \|v\|_{H^1}^2+\frac14 \m_2(u)-C.
\end{align*}
Thus, by the Banach-Alaoglu, Arzel\`a-Ascoli and Prokhorov theorems, a minimizing sequence $(u_n,v_n)_{n\in\N}$ for $\calH_\tau(\cdot|\tilde u,\tilde v)$ converges -- at least on a subsequence -- to some limit $(u,v)\in \prb(\R)\times H^1(\R)$ with $\int_\R u\log u\dd x<\infty$: $v_n\rightharpoonup v$ in $H^1(\R)$, $v_n\to v$ locally uniformly in $\R$ and $u_n\rightharpoonup u$ narrowly in $\prbb(\R)$. With respect to these convergences $\calH_\tau(\cdot|\tilde u,\tilde v)$ is lower semicontinuous, which is clear except for the term $\int_\R u_nv_n\dd x$. We employ a truncation argument similar as in \cite{zinsl2013} to prove l.s.c. for this remaining term, and consequently obtain the minimizing property for $(u,v)$. It remains to prove the additional regularity estimate \eqref{eq:addreg}. We investigate the dissipation of $\calH$ along the (auxiliary) 0-flow $(\au^s,\av^s)_{s\ge 0}$ w.r.t. $\dist$ generated by the $0$-geodesically convex functional
$$\ent(u,v):=\int_\R\left[u\log u+\frac12 v_x^2+\frac{\kappa}{2}v^2\right]\dd x$$
on $\theX$. Elementary calculations yield, since we have $\au^s_s=\au^s_{xx}$, $\av^s_s=\av^s_{xx}-\kappa\av^s$:
\begin{align*}
\frac{\dd}{\dd s}\calH(\au^s(u),\av^s(v))&\le \int_\R\left[-4(\sqrt{\au^s})_x^2+\|W_{xx}\|_{L^\infty}-\frac12(\av^s_{xx}-\kappa \av^s)^2+\frac52 \chi^2(\au^s)^2+\frac{\kappa^2}{2}(\av^s)^2\right]\dd x.
\end{align*}
Using the Sobolev inequality 
\begin{align}\label{eq:sobL4}
\|\eta\|_{L^4}\le C\|\eta\|_{H^1}^{1/4}\|\eta\|_{L^2}^{3/4},
\end{align}
we eventually arrive at
\begin{align}
\label{eq:dspHE}
\frac{\dd}{\dd s}\calH(\au^s(u),\av^s(v))&\le -2\|(\sqrt{\au^s})_x\|_{L^2}^2-\frac12 \|\av^s_{xx}-\kappa\av^s\|_{L^2}^2+\frac{\kappa^2}{2}\|v\|_{L^2}^2+C_2.
\end{align}
Finally, we use the \emph{flow interchange lemma} \cite[Thm. 3.2]{matthes2009} to obtain
$\ent(u,v)+\tau\dff^\ent \calH(u,v) \le \ent(\tilde u,\tilde v),$
which yields \eqref{eq:addreg} in combination with \eqref{eq:dspHE} and lower semicontinuity as $s\searrow 0$.\qed

Proceeding as in \cite{zinsl2013, zinsl2014}, we end up with a weak solution $(u,v)$ to \eqref{eq:pde_u}--\eqref{eq:ic} with the properties
$$v\in L^\infty([0,T];L^2(\R)),~v_x\in L^\infty([0,T];L^2(\R)),~v_t\in L^2([0,T];L^2(\R)).$$
We immediately deduce that $v\in L^\infty([0,T]\times\R)$. We now show that $v$ is continuous in both arguments. In fact, for all bounded intervals $I\subset\R$, $v$ belongs to the \emph{anisotropic} Sobolev space $W^{1,\matr{P}}([0,T]\times I)$ with
$\matr{P}=\left(\begin{smallmatrix}\frac12 & 0 \\ \frac12 & \frac12\end{smallmatrix}\right),$
the spectral radius of which is less than $1$. Since in this case $W^{1,\matr{P}}([0,T]\times I)\Subset C^0([0,T]\times\overline{I})$, the claim follows (for details on anisotropic spaces, see e.g. \cite{besov1978, krejci2011}).

\section{Sketch of proof for Theorem \ref{thm:conv}}

It is easily shown that the additional assumption of $\lambda_0$-convexity of the confinement $W$ yields boundedness from below of the entropy $\calH$. We thus obtain $(\equ,\eqv)\in (\prb\cap L^\infty)(\R)\times H^2(\R)$ as minimizer of $\calH$ similar as in \cite{zinsl2014}. Uniqueness is proved by showing strict convexity of $\calH$ as a functional on $L^2(\R)\times L^2(\R)$, which requires small coupling strength $\eps>0$.

Using the properties of $(\equ,\eqv)$, we observe that the entropy can be decomposed as follows into a convex part $\calL$ (see Proposition \ref{prop:propL} below) and a non-convex, but controllable part $\eps\calL_*$:
\begin{align}
\calH(u,v)-\calH(\equ,\eqv)&=\calL(u,v)+\eps \calL_*(u,v),\label{eq:decomp}
\end{align}
where $\calL(u,v)=\calL_u(u)+\calL_v(v),$
\begin{align*}
\calL_u(u)&:=\int_\R\left[u\log u-\equ\log\equ+W^\eps(u-\equ)\right]\dd x,~\text{with}~W^\eps:=W-\eps\eqv,\\
\calL_v(v)&:=\frac12 \|(v-\eqv)_x\|_{L^2}^2+\frac{\kappa}{2} \|v-\eqv\|_{L^2}^2,\quad
\calL_*(u,v):=-\int_\R(u-\equ)(v-\eqv)\dd x.
\end{align*}

\begin{e-proposition}[Properties of $\calL$]\label{prop:propL}
Let $\eps$ be sufficiently small. Then, the following statements hold:
\begin{enumerate}[(a)]
\item There exists $M_1>0$ such that the \emph{perturbed potential} $W^\eps$ is $\lambda_\eps$-convex, where $\lambda_\eps:=\lambda_0-M_1\eps>0$.
\item The functional $\calL_u$ is $\lambda_\eps$-geodesically convex on $(\prb(\R),\W_2)$ and
\begin{align*}
\frac{\lambda_\eps}{2}\W_2^2(u,\equ)&\le \calL_u(u)\le \frac1{2\lambda_\eps}\int_\R u((\log u+W^\eps)_x)^2\dd x.
\end{align*}
\item The functional $\calL_v$ is $\kappa$-geodesically convex on $L^2(\R)$ and
\begin{align*}
\frac{\kappa}{2}\|v-\eqv\|_{L^2}^2\le \calL_v(v)\le \frac1{2\kappa}\|(v-\eqv)_{xx}-\kappa(v-\eqv)\|_{L^2}^2.
\end{align*}
\item There exists $M_2>0$ such that $\calL(u,v)\le (1+M_2\eps)(\calH(u,v)-\calH(\equ,\eqv))$. 
\end{enumerate}
\end{e-proposition}

Actually, in one spatial dimension, the proof of part (a) simplifies dramatically compared to \cite{zinsl2014}, since
\begin{align*}
W^\eps_{xx}=W_{xx}-\eps\eqv_{xx}=W_{xx}-\eps(\kappa\eqv-\eps\equ)\ge \lambda_0-\eps\kappa\|\eqv\|_{L^\infty}\ge \lambda_0-\eps\tilde C(\calH(\equ,\eqv)+1),
\end{align*}
for some constant $\tilde C>0$. The proof of part (d) is mainly a consequence of the \emph{Csisz\'ar-Kullback inequality} (cf. \cite{carrillo2001})
%\begin{align*}
$\|u-\equ\|_{L^1}^2\le C\calL_u(u).$
%\end{align*}
In fact, more is true in one dimension: There exists another constant $C'>0$ such that
\begin{align}
\label{eq:L2fish}
\|u-\equ\|_{L^2}^2\le C' \int_\R u((\log u+W^\eps)_x)^2\dd x.
\end{align}
The idea of proof of \eqref{eq:L2fish} is as follows: We distinguish the cases where the integral on the r.h.s. in \eqref{eq:L2fish} is small or large, respectively. In the former case, we can deduce from that an $L^\infty$ bound on $u$ leading to the desired result using the Taylor expansion of the integrand in $\calL_u$ at $u_\infty(x)$. The latter case can be treated by a suitable Sobolev interpolation in one spatial dimension.

We now prove the central estimate leading to Theorem \ref{thm:conv}:
\begin{e-proposition}[Exponential estimate for $\calL$]
Let $(u_\tau^n,v_\tau^n)_{n\in \N}$ be a family of time-discrete approximations obtained by \eqref{eq:jko} which converges to a weak solution $(u,v)$ as $\tau\searrow 0$. Then, there exist $\bar\eps>0$ and $L>0$ such that for all $\eps\in (0,\bar\eps)$ and $n\in\N$, one has
\begin{align}
\label{eq:EEL}
\calL(u_\tau^n,v_\tau^n)&\le (1+M_2\eps)(\calH(u^0,v^0)-\calH(\equ,\eqv))(1+2\Lambda_\eps\tau)^{-n},
\end{align}
with $\Lambda_\eps:=\min(\lambda_0,\kappa)-L\eps>0$.
\end{e-proposition}
Once proven, this result yields exponential convergence of $\calL(u(t),v(t))$ to zero for $t\to\infty$ after passage to the continuous-time limit $\tau\searrow 0$. From this, Theorem \ref{thm:conv} clearly follows.\\
\textit{Proof. }We investigate the dissipation of $\calH$ along the (auxiliary) $\min(\lambda_\eps,\kappa)$-flow $(\au^s,\av^s)_{s\ge 0}$ of the $\min(\lambda_\eps,\kappa)$-geodesically convex functional $\calL$ on $\theX$, which is associated to the evolution system
\begin{align*}
\au^s_s&=(\au^s_{x}+\au^sW^\eps_x)_x,\qquad
\av^s_s=(\av^s-\eqv)_{xx}-\kappa(\av^s-\eqv).
\end{align*}
First, by elementary calculations, we obtain using the decomposition \eqref{eq:decomp}
\begin{align*}
\frac{\dd}{\dd s}\calH(\au^s(u),\av^s(v))&\le \left(\frac{\eps}{2}-1\right)\int_\R \au^s((\log \au^s+W^\eps)_x)^2\dd x+\frac{\eps}{2}\int_\R \au^s(\av^s-\eqv)_x^2\dd x\\&+\frac{\eps}{2}\|\au^s-\equ\|_{L^2}^2+\left(\frac{\eps}{2}-1\right)\|(\av^s-\eqv)_{xx}-\kappa(\av^s-\eqv)\|_{L^2}^2.
\end{align*}
The third term can be controlled by the first one using \eqref{eq:L2fish}, whereas the second term is to be controlled by the fourth term using the inequality $\|\eta\mu_x^2\|_{L^1}\le C\|\eta\|_{L^1}\|\mu\|_{H^2}^2$ which is valid in one spatial dimension. Taking into account the properties of $\calL$ from Proposition \ref{prop:propL}, we end up with
\begin{align}
-\frac{\dd}{\dd s}\calH(\au^s(u),\av^s(v))&\ge 2(1-\eps M)\min(\lambda_\eps,\kappa)\calL(\au^s(u),\av^s(v)),\label{eq:dissHL}
\end{align}
for some constant $M>0$ if $\eps$ is sufficiently small. The application of the flow interchange lemma \cite[Thm. 3.2]{matthes2009} eventually yields with \eqref{eq:dissHL}:
%\begin{align*}
$[1+2\tau(1-\eps M)\min(\lambda_\eps,\kappa)]\calL(u_\tau^n,v_\tau^n)\le \calL(u_\tau^{n-1},v_\tau^{n-1}).$
%\end{align*}
By iteration of this estimate and Proposition \ref{prop:propL}(d), the desired estimate \eqref{eq:EEL} follows.
\qed

\def\cprime{$'$}


\begin{thebibliography}{00}

\bibitem{savare2008}
L.~Ambrosio, N.~Gigli, and G.~Savar{\'e}.
\newblock {\em Gradient Flows: In Metric Spaces and in the Space of Probability
  Measures}.
\newblock Lectures in Mathematics. Birkh{\"a}user, 2008.

\bibitem{besov1978}
O.~V. Besov, V.~P. Il{\cprime}in, and S.~M. Nikol{\cprime}ski{\u\i}.
\newblock {\em Integral representations of functions and imbedding theorems.
  {V}ol. {I}}.
\newblock V. H. Winston \& Sons, Washington, D.C.; Halsted Press , New
  York-Toronto, Ont.-London, 1978.

\bibitem{blanchet2013r}
A.~Blanchet.
\newblock A gradient flow approach to the {K}eller-{S}egel systems, 2013.
\newblock Preprint.

\bibitem{blanchet2013}
A.~Blanchet and P.~Lauren{\c{c}}ot.
\newblock The parabolic-parabolic {K}eller-{S}egel system with critical
  diffusion as a gradient flow in {$\mathbb{R}^d,\ d\ge3$}.
\newblock {\em Comm. Partial Differential Equations}, 38(4):658--686, 2013.

\bibitem{blanchet2014}
A.~Blanchet, J.~A. Carrillo, D.~Kinderlehrer, M.~Kowalczyk, P.~Lauren{\c{c}}ot,
  and S.~Lisini.
\newblock A hybrid variational principle for the {K}eller-{S}egel system in
  $\mathbb{R}^2$, 2014.
\newblock Preprint. arXiv:1407.5562.

\bibitem{carrillo2001}
J.~A. Carrillo, A.~J\"ungel, P.~A. Markowich, G.~Toscani, and A.~Unterreiter.
\newblock Entropy dissipation methods for degenerate parabolic problems and
  generalized {S}obolev inequalities.
\newblock {\em Monatshefte für Mathematik}, 133:1--82, 2001.

\bibitem{hillen2004}
T.~Hillen and A.~Potapov.
\newblock The one-dimensional chemotaxis model: global existence and asymptotic
  profile.
\newblock {\em Math. Methods Appl. Sci.}, 27(15):1783--1801, 2004.

\bibitem{horstmann2003}
D.~Horstmann.
\newblock From 1970 until present: the {K}eller-{S}egel model in chemotaxis and
  its consequences. {I}.
\newblock {\em Jahresber. Deutsch. Math.-Verein.}, 105(3):103--165, 2003.

\bibitem{jko1998}
R.~Jordan, D.~Kinderlehrer, and F.~Otto.
\newblock The variational formulation of the {F}okker-{P}lanck equation.
\newblock {\em SIAM J. Math. Anal.}, 29(1):1--17, 1998.

\bibitem{kinderlehrer2009}
D.~Kinderlehrer and M.~Kowalczyk.
\newblock The {J}anossy effect and hybrid variational principles.
\newblock {\em Discrete Contin. Dyn. Syst. Ser. B}, 11(1):153--176, 2009.

\bibitem{krejci2011}
P.~Krej{\v{c}}{\'{\i}} and L.~Panizzi.
\newblock Regularity and uniqueness in quasilinear parabolic systems.
\newblock {\em Appl. Math.}, 56(4):341--370, 2011.

\bibitem{laurencot2011}
P.~Lauren{\c{c}}ot and B.-V. Matioc.
\newblock A gradient flow approach to a thin film approximation of the {M}uskat
  problem.
\newblock {\em Calc. Var. Partial Differential Equations}, 47(1-2):319--341,
  2013.

\bibitem{matthes2009}
D.~Matthes, R.~J. McCann, and G.~Savar{\'e}.
\newblock {A Family of Nonlinear Fourth Order Equations of Gradient Flow Type}.
\newblock {\em {Communications in Partial Differential Equations}},
  {34}({11}):{1352--1397}, {2009}.

\bibitem{mimura2012}
Y.~Mimura.
\newblock The variational formulation of the fully parabolic {K}eller-{S}egel
  system with degenerate diffusion, 2012.
\newblock Preprint.

\bibitem{osaki2001}
K.~Osaki and A.~Yagi.
\newblock Finite dimensional attractor for one-dimensional {K}eller-{S}egel
  equations.
\newblock {\em Funkcial. Ekvac.}, 44(3):441--469, 2001.

\bibitem{zinsl2013}
J.~Zinsl.
\newblock Existence of solutions for a nonlinear system of parabolic equations
  with gradient flow structure.
\newblock {\em Monatsh. Math.}, 174(4):653--679, 2014.

\bibitem{zinsl2014}
J.~Zinsl and D.~Matthes.
\newblock Exponential convergence to equilibrium in a coupled gradient flow
  system modelling chemotaxis, 2014.
\newblock Preprint. arXiv:1310.3977.

\end{thebibliography}
\end{document}